\documentclass[11pt,twoside]{article}
\usepackage{amsthm}
\usepackage{amsmath}
\usepackage{amssymb}
\usepackage{array}
\usepackage{amsfonts,amscd}
\usepackage{exscale}
\usepackage{eucal}
\usepackage{latexsym,amsmath}
\usepackage{indentfirst}
\numberwithin{equation}{section}
\newtheorem{Prop}{\bf Proposition}[section]
\newtheorem{Cor}{\bf Corollary}[section]

\newtheorem{Rem}{\bf Remark}[section]

\begin{document}
\def \b{\Box}

\begin{center}
{\Large {\bf Poisson geometry of the Maxwell-Bloch top\\[0.1cm]
system and stability problem}}
\end{center}

\begin{center}
{\bf Mihai IVAN}
\end{center}

\setcounter{page}{1}

\pagestyle{myheadings}

{\small {\bf Abstract}. Dynamics of Maxwell-Bloch top system, that
includes Maxwell-Bloch and Lorenz-Hamilton equations as particular
cases, is studied in the framework Poisson geometry. Constants of
motion as well as the relation of solution to that of pendulum are
presented. Equilibrium states are determined and, their complete
stability analysis are performed. Results are applied to an
optimal control problem on  the Lie group $G_{4}$.
}

\section{Introduction}
\smallskip
\indent  The Hamilton-Poisson systems appear naturally in many
areas of physical science and engineering including theoretical
mechanics of fluids, spatial dynamics and many others \cite{arno,
mara, gunu}. A remarkable class of Hamilton-Poisson systems is
formed by a family of differential equations on ${\bf R}^{3}$
which depend by a triple of real parameters, called the
Maxwell-Bloch top system. For certain values of these parameters
various integrable systems, such as the real-valued Maxwell-Boch
equations \cite{mara}, Lorenz-Hamilton system \cite{ivan}, etc.,
are obtained. We shall show that the solution of optimal control
problem for left invariant systems on certain matrix Lie groups
leads to systems of differential equations belonging to the family
of Maxwell-Bloch top.

 This paper is structured as follows. In Section 2, we introduce the
Maxwell-Bloch top system $(2.1)$ and some dynamical properties of
it are established. Also, we show the relation between solution of
the Maxwell-Bloch top and that of a pendulum. In Section 3, we
investigate Maxwell-Bloch top system in terms of Poisson geometry.
Section 4 is dedicated to study of Lyapunov stability for
equilibrium states of Maxwell-Bloch top system. In Section 5, we
apply results of Sections 2-4 for an optimal control problem of a
particular drift-free left invariant system on the special
nilpotent four-dimensional Lie group $G_{4}.$
\smallskip

\section{Dynamical properties of the Maxwell-Bloch top system}
\smallskip
Consider the following  family of differential equations of
Maxwell-Bloch type on ${\bf R}^{3}$:\\[-0.3cm]
\begin{equation}
\dot{x}_{1}(t) =  b_{1}
 x_{2}(t),~~
  \dot{x}_{2}(t) = b_{2}
  x_{1}(t)x_{3}(t),~~
 \dot{x}_{3}(t)  =  b_{3} x_{1}(t)x_{2}(t),
 \label{(2.1)}\\[-0.1cm]
 \end{equation}
 where $ \dot{x}_{i} = d x_{i}(t)/ dt,~ i=1, 2, 3 $ and  $ b_{1}, b_{2}, b_{3} \in{\bf R} $ are
 parameters such that $
 b_{1} b_{2} b_{3}\neq 0 $ and $ t $ is the time.
 We will refer to the dynamical system $(2.1)$ as the  {\it Maxwell-Bloch
top system} and denote the vector of parameters by $ b =
 ( b_{1}, b_{2}, b_{3})$.

If in $(2.1),$ we take  $ b = (1,  1, -1) $, then  we obtain the
{\it three-dimensional real-valued
 Maxwell-Bloch equations} \cite{daho}, given by\\[-0.3cm]
\begin{equation}
\dot{x}_{1} =
 x_{2},~~~ \dot{x}_{2} =
 x_{1}x_{3},~~~\dot{x}_{3} =- x_{1}x_{2}. \label{(2.2)}\\[-0.1cm]
\end{equation}

Also, for  $ b = (1 / 2,  -1, 1) $  we obtain the {\it Lorenz-Hamilton system} \cite{gunu,ivan}, given by\\[-0.2cm]
\begin{equation}
\dot{x}_{1} =
 \dfrac{1}{2}x_{2},~~~ \dot{x}_{2} =-
 x_{1}x_{3},~~~\dot{x}_{3} =  x_{1}x_{2}. \label{(2.3)}\\[-0.1cm]
\end{equation}
\begin{Prop}
 The functions $ H^{b}, C^{b} \in C^{\infty}({\bf R}^{3}, {\bf R}) $  given by:\\[-0.2cm]
\begin{equation}
 H^{b}(x_{1},x_{2}, x_{3}) =
\displaystyle\frac{b_{1}}{2}\left(x_{2}^{2}-
\displaystyle\frac{b_{2}}{b_{3}} x_{3}^{2}\right)~~~\hbox{and}~~~
C^{b}(x_{1}, x_{2}, x_{3}) = -\displaystyle\frac{b_{3}}{2 b_{1}}
x_{1}^{2} + x_{3} \label{(2.4)}
\end{equation}
are constants of the motion (first integrals) for the dynamics
(2.1).
\end{Prop}
{\it Proof.} Indeed,\\[0.1cm]
 $\displaystyle\frac{d H^{b}}{dt} = b_{1}(
x_{2} \dot{x}_{2}- \dfrac{b_{2}}{b_{3}} x_{3}\dot{x}_{3})= b_{1}(
- x_{2}(b_{2}x_{1}x_{3}) + \dfrac{b_{2}}{b_{3}}
x_{3}(b_{3} x_{1}x_{2}))= 0~$ and\\[0.1cm]
$\displaystyle\frac{dC^{b}}{dt} = -\dfrac{b_{3}}{b_{1}} x_{1}
\dot{x}_{1}+\dot{x}_{3}= -\dfrac{b_{3}}{b_{1}}x_{1} (b_{1}x_{2}) +
b_{3}x_{1}x_{2}= 0 .~~~\hfill\Box$

\begin{Rem} {\rm From Proposition $2.1$ it follows that} the trajectories of the dynamics (2.1) in
the phase space ${\bf R}^{3}$ are the intersections of the
surfaces:\\[-0.3cm]
\[
 b_{1} x_{2}^{2} - \displaystyle\frac{b_{1}b_{2}}{b_{3}}
 x_{3}^{2}= 2H^{b},~~~-\displaystyle\frac{b_{3}}{b_{1}}
x_{1}^{2} + 2 x_{3}= 2C^{b},\\[-0.1cm]
\]
where $~H^{b}=constant~$ and $~C^{b}=constant .$  \hfill$\Box$
\end{Rem}

\markboth{M. Ivan}{Poisson geometry of the Maxwell-Bloch top
system and stability problem}

Using the fact that  $H^{b}$ given by $(2.4)$ is a first integral
(see Proposition $2.1$) one easily prove that the Maxwell-Bloch
top system has the following first integral:
\begin{equation}
H_{0}^{b}(x)= \displaystyle\frac{1}{2}\left (x_{2}^{2} -
\displaystyle\frac{b_{2}}{b_{3}} x_{3}^{2}\right ). \label{(2.5)}
\end{equation}

We shall prove that in certain restrictions on  $b_{i}$, the
motion of Maxwell-Bloch top system reduces to motion on the
surface described by the conservation law $(2.5)$.
\begin{Prop}
We assume that $b_{2} b_{3} < 0.$ The solution of the
Maxwell-Bloch top system $(2.1)$ restricted to the constant level
surface defined by:\\[-0.3cm]
\begin{equation}
x_{2}^{2} - \displaystyle\frac{b_{2}}{b_{3}} x_{3}^{2}= 2 H =
constant, ~~~ H = H_{0}^{b} > 0 \label{(2.6)}\\[-0.1cm]
\end{equation}
is\\[-0.5cm]
\begin{equation}
\left\{ \begin{array}{lcll}
 x_{1}(t)& =&\dfrac{\gamma}{b_{3}}\cdot \dot{\theta}(t)&~~\hbox{with}~~\gamma=\sqrt{-\dfrac{b_{3}}{b_{2}}}\\[0.3cm]

x_{2}(t) &=&
\sqrt{2H}\cdot \cos\theta(t)& \\[0.1cm]

x_{3}(t) &=& \gamma \sqrt{2H}\cdot \sin \theta(t)&\\
\end{array}\right. \label{(2.7)}\\[-0.1cm]
\end{equation}
where $\theta(t)$ is a solution of the pendulum equation:\\[-0.3cm]
\begin{equation}
\ddot{\theta}(t)= \dfrac{b_{1}b_{3}}{\gamma}\sqrt{2H} \cdot \cos
\theta(t).\label{(2.8)}\\[-0.1cm]
\end{equation}
\end{Prop}
{\it Proof.} Denote $\gamma=\sqrt{-b_{3} / b_{2}} >0.$  By a
direct computation, it is easy to see that

$(i)~~~~~~~ x_{2}(t) =\sqrt{2H}\cdot \cos\theta(t),~~~ x_{3}(t)
=\gamma \sqrt{2H}\cdot \sin \theta(t)$\\[0.1cm]
are solutions of the equation  $(2.6)$. By deriving of the second
relation of $(i)$ with respect to  $t$, we have
 $~\dot{x}_{3}(t) = \gamma \sqrt{2H}\cdot
\cos\theta(t)\cdot \dot{\theta}(t)~$ and using the first relation
of $(i)$, we obtain:\\[-0.2cm]

$(ii)~~~~~~~ \dot{x}_{3}(t) =\gamma \cdot x_{2}(t)\cdot
\dot{\theta}(t).$\\[-0.2cm]

From $(ii)$  and $\dot{x}_{3}(t) = b_{3} x_{1}(t) x_{2}(t)$, we
deduce
 $~x_{1}(t) = \displaystyle\frac{\gamma}{b_{3}}\cdot
\dot{\theta}(t).$ Therefore, the relations $(2.7)$ are verified.
From the last equality follows:\\[-0.2cm]

$(iii)~~~~~~~\dot{\theta}(t) = \dfrac{b_{3}}{\gamma}\cdot
x_{1}(t).$\\[-0.2cm]

Differentiating again the relation $(iii)$ and using the first
equation from $(2.1)$ and $(i)$, it follows $~\ddot{\theta}(t)=
\dfrac{b_{3}}{\gamma}\cdot \dot{x}_{1}(t)=
\dfrac{b_{1}b_{3}}{\gamma} \sqrt{2 H}\cdot \cos \theta(t),~$ i.e.
$(2.8)$ holds.\hfill$\Box$

\begin{Cor}
The solution of the Lorenz-Hamilton system $(2.3)$, restricted to
the constant level surface defined by:\\[-0.3cm]
\[
x_{2}^{2} + x_{3}^{2}= 2 H = constant, ~~~ H > 0\\[-0.1cm]
\]
is\\[-0.5cm]
\begin{equation}
 x_{1}(t) = \dot{\theta}(t),~~~
x_{2}(t) = \sqrt{2H}\cdot \cos\theta(t),~~~ x_{3}(t) =
\sqrt{2H}\cdot \sin \theta(t),\label{(2.9)}
\end{equation}
where $\theta(t)$ is a solution of the pendulum
equation:\\[-0.2cm]
\begin{equation}
\ddot{\theta}(t)= \dfrac{1}{2}\sqrt{2H} \cdot \cos
\theta(t).\label{(2.10)}\\[-0.2cm]
\end{equation}
\end{Cor}
{\it Proof.} In Proposition 2.2 we take $b_{1}=1 / 2,~ b_{2}=-1,~
b_{3} = 1$ and we obtain the required result.\hfill$\Box$
\smallskip

\section{Realizations Hamilton-Poisson for the Maxwell-Bloch top system}

\smallskip
 For definitions and results on Poisson geometry and
Hamiltonian dynamics see \cite{mara, gunu}.

\begin{Prop}
$(i)~$ The Maxwell-Bloch top system (2.1) is a Hamilton-Poisson
system with the phase space ${\bf R}^{3}$, the Hamiltonian $ H^{b}
$ given by $(2.4)$ and with respect the Poisson structure
$\{\cdot, \cdot\} $ given
by\\[-0.2cm]
\begin{equation}
\{f,g\} = \det\left(\begin{array}{ccc} -\dfrac{b_{3}}{b_{1}}x_{1}
& 0 & 1\\ \\[-0.2cm]
 \displaystyle\frac{\partial f}{\partial x_{1}}&
\displaystyle\frac{\partial f}{\partial
x_{2}}&\displaystyle\frac{\partial f}{\partial x_{3}}\\ \\[-0.2cm]
\displaystyle\frac{\partial g}{\partial x_{1}}&
\displaystyle\frac{\partial g}{\partial
x_{2}}&\displaystyle\frac{\partial g}{\partial x_{3}}\\
\end{array}\right ),~~\hbox{for all}~~ f,g\in C^{\infty}({\bf
R}^{3}).\label{(3.1)}
\end{equation}

(ii)$~$ The function $ C^{b}$ given by $(2.4)$ is a Casimir of the
configuration $ ( {\bf R}^{3}, \{\cdot, \cdot\})$.
\end{Prop}

{\it Proof.} $(i)~$ It is easy to observe that $\{f,g\}= \nabla
C^{b}\cdot (\nabla f\times \nabla g).$ Then $\{\cdot, \cdot\}$ is
a bracket operation on ${\bf R}^{3}.$

The system $(2.1)$ is a Hamilton-Poisson system, since
$~\dot{x}_{i} = \{x_{i}, H^{b}\},~i=1, 2, 3.$\\[0.1cm]
 Indeed, for
example\\[0.1cm]
 $\{x_{1}, H^{b}\} =
\left|\begin{array}{ccc} -\dfrac{b_{3}}{b_{1}}x_{1} & 0 & 1\\
\\[-0.2cm]
1 & 0&0\\ \\[-0.2cm]
0 &  b_{1} x_{2} & -\frac{b_{1}b_{2}}{b_{3}} x_{3}\\
\end{array}\right | = b_{1} x_{2} = \dot{x}_{1}. $\\

 $(ii)~$ The function $C^{b}\in C^{\infty}({\bf R}^{3}, \bf {\bf R})$ is  a Casimir, since  $\{C^{b}, f\}=0$ for every $ f\in C^{\infty}({\bf
R}^{3}, {\bf R}).$ We have $~\{C^{b},f\} =\nabla C^{b}\cdot
(\nabla C^{b}\times \nabla f)= 0.$ \hfill$\Box$

We can easily prove that the Poisson structure $\{\cdot,\cdot\}$
given by $(3.1)$ is in fact generated by the skew-symmetric matrix
\begin{equation}
P^{b}(x_{1}, x_{2}, x_{3}) = \left ( \begin{array}{ccc} 0 & 1 & 0\\ \\[-0.2cm]
-1 & 0 & -\dfrac{b_{3}}{b_{1}}x_{1}\\ \\[-0.2cm]
0 &  \dfrac{b_{3}}{b_{1}}x_{1}& 0\\
\end{array}\right ).\label{(3.2)}
\end{equation}

The Maxwell-Bloch top system $(2.1)$ can be expressed in the
matrix form:
\begin{equation}
\dot{X} = P^{b}(x)\cdot \nabla H(x), \label{(3.3)}
\end{equation}
where $~x=(x_{1}, x_{2}, x_{3}) $ and $~\dot{X}
=\left(\begin{array}{cccc}
 \dot{x}_{1} & \dot{x}_{2}& \dot{x}_{3}\\
\end{array}\right )^{T}.$

 Define the functions $C_{\alpha \beta}^{b}, H_{\gamma\delta}^{b}\in
C^{\infty}({\bf R}^{3}, {\bf R})$ be given by:
\begin{equation}
C_{\alpha\beta}^{b} = \alpha C^{b} + \beta H^{b},~~~
H_{\gamma\delta}^{b} =\gamma C^{b} + \delta H^{b},~~\alpha, \beta,
\gamma,\delta\in {\bf R} ~~~~~\hbox{that is}\label{(3.4)}
\end{equation}
\begin{equation}
\left\{ \begin{array}{ccc}
C_{\alpha\beta}^{b}(x_{1}, x_{2}, x_{3}) &=& -\displaystyle\frac{\alpha b_{3}}{2 b_{1}} x_{1}^{2} + \displaystyle\frac{\beta b_{1}}{2}x_{2}^{2}
+ \alpha x_{3} - \displaystyle\frac{\beta b_{1}b_{2}}{2 b_{3}} x_{3}^{2} \\[0.4cm]
H_{\gamma\delta}^{b}(x_{1}, x_{2}, x_{3}) &=&
-\displaystyle\frac{\gamma b_{3}}{2 b_{1}} x_{1}^{2} +
\displaystyle\frac{\delta b_{1}}{2}x_{2}^{2}
+ \gamma x_{3} - \displaystyle\frac{\delta b_{1}b_{2}}{2 b_{3}} x_{3}^{2} \\[-0.2cm]
\end{array}\right.\label{(3.5)}
\end{equation}
\begin{Prop}
$(i)$ The Maxwell-Bloch top system $(2.1)$ admits a family of
Hamilton-Poisson realizations parametrized by the group $SL(2;{\bf
R})$. More precisely, $({\bf R}^{3}, \{\cdot,
\cdot\}_{\alpha\beta}^{b}, H_{\gamma\delta}^{b})$ is a
Hamilton-Poisson realization of the system $(2.1)$, where
$\{\cdot, \cdot\}_{\alpha\beta}$ is given by:
\begin{equation}
\{f,g\}_{\alpha\beta}^{b} = \det\left(\begin{array}{ccc}
-\dfrac{\alpha b_{3}}{b_{1}}x_{1} & \beta b_{1} x_{2} & \alpha- \dfrac{\beta b_{1}b_{2}}{b_{3}}x_{3}\\
\\[-0.2cm]
\displaystyle\frac{\partial f}{\partial x_{1}}&
\displaystyle\frac{\partial f}{\partial
x_{2}}&\displaystyle\frac{\partial f}{\partial x_{3}}\\ \\[-0.2cm]
\displaystyle\frac{\partial g}{\partial x_{1}}&
\displaystyle\frac{\partial g}{\partial
x_{2}}&\displaystyle\frac{\partial g}{\partial x_{3}}\\
\end{array}\right ), ~ \forall f,g\in C^{\infty}({\bf
R}^{3},{\bf R}),\label{(3.6)}\\[0.2cm]
\end{equation}
 the Hamiltonian $H_{\gamma\delta}^{b} $ is given by $(3.5)$ and the matrix
 $\left(\begin{array}
   {cc}
   \alpha &\beta\\
   \gamma &\delta
   \end{array}\right)\in SL(2;{\bf R})$.

$(ii)~ C_{\alpha\beta}^{b}$ given by $(3.5)$ is a Casimir of the
configuration $ ( {\bf R}^{3}, \{\cdot, \cdot
\}_{\alpha\beta}^{b})$.
\end{Prop}

{\it Proof.} $(i)~$ We have $ \displaystyle\frac{\partial
H_{\gamma\delta}^{b}}{\partial x_{1}}=
 \dfrac{\gamma b_{3}}{b_{1}}x_{1},~\displaystyle\frac{\partial H_{\gamma\delta}^{b}}{\partial
x_{2}}=\delta b_{1} x_{2},~\displaystyle\frac{\partial
H_{\gamma\delta}^{b}}{\partial x_{3}}=\gamma -\dfrac{\delta b_{1}
b_{2}}{b_{3}} x_{3}.~$ Then:
\[
 \{x_{1} , H_{\gamma\delta}^{b}\}_{\alpha\beta}^{b}= \det \left (\begin{array}{ccc}
-\dfrac{\alpha b_{3}}{b_{1}}x_{1} & \beta b_{1} x_{2} & \alpha- \dfrac{\beta b_{1}b_{2}}{b_{3}}x_{3}\\ \\[-0.2cm]
1 & 0 & 0\\ \\[-0.2cm]
 \dfrac{\gamma b_{3}}{b_{1}}x_{1}& \delta b_{1} x_{2}
& \gamma - \dfrac{\delta
b_{1} b_{2}}{b_{3}} x_{3}\\
 \end{array}\right )= (\alpha\delta - \beta\gamma) b_{1} x_{2} = \dot{x}_{1}.
\]

Similarly, we have $~\{x_{2},
H_{\gamma\delta}^{b}\}_{\alpha\beta}^{b} = b_{2} x_{1} x_{3} =
\dot{x}_{2}~$ and $~\{x_{3},
H_{\gamma\delta}^{b}\}_{\alpha\beta}^{b} = b_{3} x_{1} x_{2} =
\dot{x}_{3}.$ Therefore, one obtains the required result.

$(ii)~$ It is easy to see that
$\{C_{\alpha\beta}^{b},f\}_{\alpha\beta}^{b} = 0, $ for all $ f
\in C^{\infty}({\bf R}^{3}, {\bf R}).$\hfill$\Box$\\

 The Poisson structure given by $(3.6)$ is generated by the matrix
\begin{equation}
P_{\alpha\beta}^{b}(x_{1},x_{2}, x_{3})= \left (\begin{array}{ccc}
0 & \alpha-\dfrac{\beta b_{1} b_{2}}{b_{3}} x_{3}& - \beta b_{1}
x_{2}\\ \\[-0.2cm]
-\alpha+\dfrac{\beta b_{1} b_{2}}{b_{3}} x_{3} & 0 & -
\dfrac{\alpha b_{3}}{b_{1}} x_{1}\\ \\
\beta b_{1}x_{2}& \dfrac{\alpha b_{3}}{b_{1}} x_{1} &0\\ \\[-0.2cm]
\end{array}\right ).\label{(3.7)}
\end{equation}
\begin{Rem}
{\rm Proposition 3.2 assures that the equations $(2.1)$ are
invariant, if  $ H^{b}$ and $C^{b}$ are replaced by linear
combinations with coefficients modulo $ SL(2,{\bf R})$. In
consequence, the trajectories of motion of the system $(2.1)$
remain unchanged.}\hfill$\Box$
\end{Rem}

Finally, we can conclude that {\it the Maxwell-Bloch top system
$(2.1)$ has the following Hamilton-Poisson realization\\[-0.2cm]
\begin{equation}
({\bf R}^{3}, P_{\alpha\beta}^{b},
H_{\gamma\delta}^{b})~~~\hbox{with
Casimir}~~~C_{\alpha\beta}^{b},\label{(3.8)}
\end{equation}
where $~P_{\alpha\beta}^{b}~$ is given by $(3.7)~$ and
$~H_{\gamma\delta}^{b},~C_{\alpha\beta}^{b}~$ are given by $(3.5)$
for all $\alpha, \beta, \gamma, \delta \in {\bf R} $ such that
$\alpha\delta-\beta\gamma =1.$}

 If in $(3.8)$ we take $ \alpha
=1, \beta = \gamma = 0 $ and $\delta = 1$, then one obtains
Proposition 3.1. More precisely, {\it $~({\bf R}^{3}, P^{b},
H^{b})~$  is a Hamilton-Poisson realization of the dynamics
$(2.1)$ with Casimir $~C^{b}$},
 since $~H_{01}^{b} =
H^{b},~C_{10}^{b} = C^{b} $ and $ P_{10}^{b} = P^{b}.$

Next proposition gives another special Hamilton-Poisson
realization of the Maxwell-Bloch top system.

\begin{Prop}
$(i)~$ The Maxwell-Bloch top  system (2.1) has the
Hamilton-Poisson realization $~({\bf R}^{3}, \bar{P}^{b},
\bar{H}^{b}),~$ where the matrix $~\bar{P}^{b}~$ is given by
\begin{equation}
\bar{P}^{b}(x_{1}, x_{2}, x_{3}) =  \left ( \begin{array}{ccc} 0 &
\displaystyle\frac{b_{1}b_{2}}{b_{3}}x_{3} & b_{1}x_{2}\\ \\[-0.2cm]
-\displaystyle\frac{b_{1}b_{2}}{b_{3}}x_{3} & 0 & 0\\ \\[-0.2cm]
-b_{1}x_{2}& 0 & 0
\end{array}\right ),\label{(3.9)}
\end{equation}
and the Hamiltonian $\bar{H}^{b}$ is given by\\[-0.2cm]
\begin{equation}
\bar{H}^{b}(x_{1}, x_{2}, x_{3}) = -\displaystyle\frac{b_{3}}{2
b_{1}} x_{1}^{2} + x_{3}. \label{(3.10)}\\[-0.1cm]
\end{equation}

 $(ii)~$ The function $ \bar{C}^{b}$ defined by\\[-0.2cm]
\begin{equation}
\bar{C}^{b}(x_{1}, x_{2}, x_{3})=
\displaystyle\frac{b_{1}}{2}\left( x_{2}^{2}-
\displaystyle\frac{b_{2}}{b_{3}} x_{3}^{2}\right).
\label{(3.11)}\\[-0.1cm]
\end{equation}
is a Casimir of the configuration $ ( {\bf R}^{3},
\{\cdot,\cdot\}_{1} )$, where $\{\cdot,\cdot\}_{1} $ is the
bracket operation whose its Poisson matrix is $\bar{P}^{b}$.
\end{Prop}

{\it Proof.} The assertions are consequences of Proposition 3.2.
For $ \beta =-1, \alpha =\delta = 0 $ and $\gamma = 1$ one obtains
$~\{\cdot,\cdot\}_{1}=\{\cdot, \cdot\}_{0,-1}^{b}, \bar{H}^{b} =
H_{10}^{b}, \bar{C}^{b}= -C_{0,-1}^{b} $ and $ \bar{P}^{b}=
P_{0,-1}^{b}.\hfill\Box$

Using Propositions 3.2 and 3.3 we obtain the following
proposition.
\begin{Prop}
The Maxwell-Bloch top system $(2.1) $ have the following two
(special) Hamilton-Poisson realizations:

$(i)~~~ ( {\bf R}^{3}, P^{b}, H^{b} ) $ with the Casimir  $
C^{b}\in C^{\infty}({\bf R}^{3}, {\bf R})$, where $P^{b}$ is given
by $(3.2)$ and $H^{b}, C^{b}$ are given by $(2.4)$;

$(i)~~~ ( {\bf R}^{3}, \bar{P}^{b}, \bar{H}^{b} ) $ with the
Casimir $ \bar{C}^{b}\in C^{\infty}({\bf R}^{3},{\bf R})$, where
$\bar{P}^{b}$ is given by $(3.9)$ and $\bar{H}^{b}, \bar{C}^{b}$
are given by $(3.10)$ and $(3.11),$ respectively.\hfill$\Box$
\end{Prop}
\begin{Rem}
{\rm We have $\bar{P}^{b} = \Pi_{(0,u,v)} $ with $u=b_{1}, v=-
b_{1} b_{2} / b_{3}$ (see the relation (2.5) in \cite{giva}),
where $~\Pi_{(0,u,v)}= \left ( \begin{array}{ccc}
 0 & -v x_{3} & u x_{2}\\
v x_{3} & 0 & 0\\
- u x_{2} & 0 & 0\\
\end{array}\right ).~$ Hence the Poisson geometry of the system (2.1) is generated by a
matrix of $se(2)-$type (here, $ se(2)$ is the Lie algebra of the
Lie group $~SE(2;{\bf R}))$.}\hfill$\Box$
\end{Rem}
\begin{Rem}
{\rm Applying Proposition 3.4 one obtains two special
Hamilton-Poisson realizations for the real Maxwell-Bloch equations
$(2.2)$ and Lorenz-Hamilton system $(2.3),$
respectively.}\hfill$\Box$
\end{Rem}

\section{ Stability problem for Maxwell-Bloch top dynamics}

A direct computation shows that the equilibrium states of the
Maxwell-Bloch top system  $(2.1)$ are the points\\[-0.2cm]
 \[
 e_{0} = (0,
0, 0),~~ e_{1}^{m} = (m, 0, 0)~~\hbox{and}~~ e_{3}^{m}= (0, 0,
m)~~\hbox{for all}~~ m \in {\bf R}^{\ast}.
\]

Let $A(x_{1}, x_{2}, x_{3})$ be the matrix of the linearisation of
the system $(2.1)$,
i.e.\\[-0.2cm]
\[
A(x_{1}, x_{2}, x_{3})= \left (\begin{array}{ccc}
  0 & b_{1}   & 0 \\
  b_{2}x_{3}  & 0 &  b_{2} x_{1}\\
   b_{3}x_{2}   &  b_{3}x_{1}  & 0 \\
\end{array}\right ).
\]

\begin{Prop}
$(i)~$ The equilibrium states  $ e_{1}^{m},~ m \in {\bf R}^{\ast}
$ are spectrally stable if $~b_{2}b_{3}< 0~$ and unstable if
$~b_{2}b_{3}> 0.$

$(ii)~$ The equilibrium states $ e_{3}^{m},~ m \in {\bf R}^{\ast}
$  are spectrally stable if $~m b_{1} b_{2}< 0~ $ and unstable if
$~m b_{1} b_{2}> 0 .$

$(iii)~$ The equilibrium state $ e_{0}=(0, 0, 0)$  is spectrally
stable.
\end{Prop}
{\it Proof.} $(i)~$ The characteristic polynomial of
 \[
A(e_{1}^{m}) =\left (\begin{array}{ccc}
  0 & b_{1} & 0\\
  0 & 0 & m b_{2}\\
  0 & m b_{3} & 0 \\
\end{array}\right )
\]  is
$~p_{A(e_{1}^{m})}(\lambda) = \det (A(e_{1}^{m})- \lambda I)= -
\lambda ( \lambda ^{2} - b_{2} b_{3}m^{2} ).$ Then the
characteristic roots of $ A(e_{1}^{m}) $ are $ \lambda_{1} = 0 $
and $ \lambda_{2,3}= \pm m \sqrt{b_{2} b_{3}},$ if
 $ b_{2} b_{3}> 0 $ and $\lambda_{2,3}= \pm
im\sqrt{-b_{2} b_{3}}, $ if $ b_{2} b_{3}< 0.$ From Lyapunov's
Theorem it follows that $ e_{1}^{m}$ is spectrally stable for $
b_{2} b_{3}< 0 $ and unstable for $ b_{2} b_{3}> 0 .$

$(ii)~$ The characteristic polynomial of
\[
A(e_{3}^{m})= \left (\begin{array}{ccc}
  0 & b_{1} & 0 \\
  m b_{2} & 0 &  0\\
  0 & 0 & 0 \\
\end{array}\right )
\]
 is $~p_{A(e_{3}^{m})}(\lambda) = - \lambda ( \lambda ^{2} - m b_{1}
b_{2} ) $ with the characteristic roots $~\lambda_{1}=0,\,
\lambda_{2,3}= \pm \sqrt{m b_{1} b_{2}}.$ Applying now similar
arguments as in the proof of the assertion $(i)$, one obtains the
required results.

$(iii)~$  It is easy to see that $ e_{0}$ is spectrally
stable.\hfill$\Box$

Let us discuss the nonlinear stability of equilibrium states of
the dynamics $(2.1)$ which are spectrally stable. Recall that an
equilibrium state $x_{e}$ is nonlinear stable if the trajectories
starting close to $x_{e}$ stay close to $x_{e}$ (i.e. a
neighborhood of $x_{e}$ must be flow invariant).

\begin{Prop}
If $ b_{2} b_{3}< 0, $ then $e_{1}^{m},$ $m\in {\bf R}^{\ast} $ is
nonlinear stable.
\end{Prop}
{\it Proof.} We suppose that $ b_{2} b_{3} <0.~$  We shall make
the proof using Lyapunov's theorem  \cite{hism}.  Let be the
function $L^{\alpha}: {\bf R}^{3}\to  {\bf R} $ given
by:\\[-0.2cm]
\[
L^{b}(x_{1}, x_{2}, x_{3})= \dfrac{1}{2}\left(
-\displaystyle\frac{b_{3}}{2b_{1}} x_{1}^{2} + x_{3} +
\displaystyle\frac{b_{3}}{2b_{1}}m^{2}\right )^{2} +
\dfrac{1}{2}\left( x_{2}^{2} - \displaystyle\frac{b_{2}}{b_{3}}
x_{3}^{2}\right).
\]

 For the function $L^{b}$ we have successively:\\[-0.2cm]

$(i)~~~ L^{b}\in C^{\infty}({\bf R}^{3}, {\bf R})~$ and $~L^{b}(m,0,0)=
0;$\\[-0.2cm]

$(ii)~~L^{b}(x_{1}, x_{2}, x_{3})> 0$, for all $ x\in {\bf R}^{3},~ x\neq e_{1}^{m},$ since $ b_{2} b_{3} <
0;$\\[-0.2cm]

$(iii)~$ The derivative of  $L^{b}$ with respect to $t$ along the
trajectories of the dynamics $(2.4)$ is zero. Indeed,
\[
\displaystyle\frac{d L^{b}}{dt}= \displaystyle\frac{\partial
L^{b}}{\partial x_{1}}\dot{x}_{1} + \displaystyle\frac{\partial
L^{b}}{\partial x_{2}}\dot{x}_{2} +\displaystyle\frac{\partial
L^{b}}{\partial x_{3}}\dot{x}_{3} = -
\displaystyle\frac{b_{3}}{b_{1}} x_{1} \left(
-\displaystyle\frac{b_{3}}{2b_{1}} x_{1}^{2} + x_{3} +
\displaystyle\frac{b_{3}}{2b_{1}}m^{2}\right)  b_{1} x_{2}+  x_{2}
( b_{2} x_{1} x_{3}) +
\]
\[
+ [ \left( -\displaystyle\frac{b_{3}}{2b_{1}} x_{1}^{2} + x_{3} +
\displaystyle\frac{b_{3}}{2b_{1}}m^{2}\right) -
\dfrac{b_{2}}{b_{3}} x_{3} ] ( b_{3} x_{1} x_{2}) =0.
\]

 Therefore $L^{b}$ is a Lyapunov function. Then via Lyapunov's theorem we obtain that $e_{1}^{m}$ is
nonlinear stable.\hfill$\Box$

\begin{Prop}
 The equilibrium state $ e_{0}$ of the dynamics $(2.1)$ is nonlinear
stable.
\end{Prop}
{\it Proof.} An easy computation shows that
  \[
   L_{0}^{b}(x_{1},
x_{2}, x_{3}) = \displaystyle\frac{1}{4}\left (x_{2}^{2} -
\displaystyle\frac{b_{2}}{b_{3}} x_{3}^{2}\right )^{2}
\]
is a Lyapunov function. The assertion is a consequence of the
Lyapunov theorem.\hfill$\Box$

\begin{Prop}
If $ m b_{1} b_{2} < 0$, then $e_{3}^{m}$ is nonlinear stable.
\end{Prop}
{\it Proof.} We shall make the proof using Arnold's energy-Casimir
method \cite{arno}. Let the function  $F_{\lambda}^{b}\in
C^{\infty}({\bf R}^{3}, {\bf R}),$ $\lambda\in {\bf R} $ given
by:\\[-0.2cm]
\[
F_{\lambda}^{b}(x_{1}, x_{2}, x_{3})= H^{b}(x_{1}, x_{2}, x_{3})-
\lambda C^{b}(x_{1}, x_{2}, x_{3}) =
\]
\[
=\displaystyle\frac{ b_{1}}{2}\left(x_{2}^{2}-
\displaystyle\frac{b_{2}}{b_{3}} x_{3}^{2}\right)- \lambda
\left(-\displaystyle\frac{b_{3}}{2 b_{1}} x_{1}^{2} +
x_{3}\right).
\]

 Then we have successively:\\[-0.2cm]

$(i)~~ \nabla F_{\lambda}^{b}(e_{3}^{m})=0~$ if and only if
$~\lambda =\lambda_{0}$, where $\lambda_{0}= -m b_{1}
b_{2}/ b_{3};$\\[-0.2cm]

$(ii)~~W:= ker~d C_{2}^{b}(e_{3}^{m}) = span \left((1, \,
0, \, 0)^{T}, (0, \, 1, \, 0)^{T}\right );$\\[-0.2cm]

$(iii)~$ For all  $ v\in W$, i.e. $ v=(\alpha, \beta, 0)^{T},~
\alpha, \beta\in {\bf R},$ we have:\\[-0.2cm]
\[
v^{T}\cdot \nabla^{2}F_{\lambda_{0}}^{b}(e_{3}^{m})\cdot v
=\dfrac{1}{b_{1}}\left ( -m b_{1} b_{2} \alpha^{2}+ b_{1}^{2}
\beta^{2}\right )
\]
and so $~\nabla^{2}F_{\lambda_{0}}^{b}(e_{3}^{m})\Big|_{W\times
W}$ is positive  definite (respectively, negative definite) if $~
b_{1}> 0 $ and $ m b_{1}b_{2} < 0~$ (respectively, $~ b_{1}<0 $
and $~m b_{1}b_{2}<0 $). Therefore via Arnold's energy-Casimir
method we conclude that $e_{3}^{m}$ is nonlinear
stable.\hfill$\Box$

\begin{Cor}
The equilibrium states  $ e_{0}$ and $e_{1}^{m}, e_{3}^{m}$ for $
m\in {\bf R}^{\ast}$ of the Lorenz-Hamilton system given by
$(2.3)$ have the following behavior:

$(i)~~~e_{0}$ and $ e_{1}^{m}$  are nonlinear stable;

$(ii)~~ e_{3}^{m} $ is nonlinear stable for $ m>0$ and  unstable
for $m<0$.
\end{Cor}
{\it Proof.} The assertions follows immediately from Propositions
4.2-4.4.\hfill$\Box$

\section{Application to study of an invariant controllable system on  $G_{4}$}

Control systems with state evolving on a matrix Lie group arise
frequently in physical problems and many others \cite{jusu, kris}.

In this section we present a drift-free left invariant
controllable system on a particular Lie group. This arise
naturally from the study of the car's dynamics for which the Lie
group $G_{4}$ represents the configuration space \cite{puta}.

We denote by $G_{4}\subset UP(4)$ the subgroup of unipotent
matrices consisting of elements $X$ of the form:\\[-0.2cm]
\[
X=\left(\begin{array}{cccc} 1 & x_{2} & x_{3} & x_{4}\\ [0.2cm]
 0&1 & x_{1} & \dfrac{x_{1}^{2}}{2}\\ [0.2cm]
0 & 0 & 1 & x_{1}\\
0 & 0 & 0 & 1\\
\end{array}\right ),~ (x_{1}, x_{2}, x_{3}, x_{4})\in {\bf R}^{4}.
\]

A basis of the Lie algebra  ${\cal G}_{4}$ associated to $G_{4}$
is $\{ A_{1}, A_{2}, A_{3}, A_{4}\}$, where:\\[-0.2cm]
\[
A_{1}=\left(\begin{array}{cccc}
0 & 0 & 0 & 0\\
0 & 0 & 1 & 0\\
0 & 0 & 0 & 1\\
0 & 0 & 0 & 0\\
\end{array}\right ),~A_{2}=\left(\begin{array}{cccc}
0 & 1 & 0 & 0\\
0 & 0 & 0 & 0\\
0 & 0 & 0 & 0\\
0 & 0 & 0 & 0\\
\end{array}\right ),~ A_{3}=[A_{2}, A_{1}],~
A_{4}=[A_{3}, A_{1}].
\]

Using the results from \cite{jusu} of controllability for
drift-free left invariant systems it follows that there exist only
four drift-free left invariant controllable systems on $G_{4}$,
namely
\begin{equation}
\dot{X}= X(A_{1} u_{1} + A_{2} u_{2} + p A_{3} u_{3} + q A_{4}
u_{4}),~~~~~X\in G_{4}\label{(5.1)}
\end{equation}
where $u_{i}\in C^{\infty}({\bf R}, {\bf R})~$ are control
 functions and $(p,q)\in \{(0,0),~(1,0),~(0,1),~(1,1)\}.$

An optimal control problem for the system $(5.1)$ with
$(p,q)=(1,0)$ has been studied in the paper \cite{poar}.

 The space of configurations for kinematics of a car is ${\bf R}^{2}\times S^{1}\times
S^{1}$, and its dynamics is described by the system of
differential equations :
\begin{equation}
\dot{x}_{1} = u_{1},~~~\dot{x}_{2} =
u_{2},~~~x_{3}=u_{1}x_{2},~~~\dot{x}_{4} = u_{1} x_{3}.
\label{(5.2)}
\end{equation}

The system $(5.2)$ can be interpreted as a drift-free left
invariant control system on $G_{4}$ \cite{jusu, puta}. Indeed, the
system $(5.2)$ can be written in the equivalent form:
\begin{equation}
\dot{X}= X( A_{1} u_{1} + A_{2}u_{2}),~~~\hbox{where}~~~X\in G_{4}
.\label{(5.3)}
\end{equation}

For the system $(5.3)$ we consider the cost function $J$ be given
by:
\begin{equation}
J(u_{1},u_{2})=\dfrac{1}{2} \int_{0}^{t_{f}} [c_{1} u_{1}^{2}(t)+
c_{2} u_{2}^{2}(t)]dt,~~~ c_{1}>0, c_{2}>0 .\label{(5.4)}
\end{equation}

Using the Krishnaprasad's theorem \cite{kris} we obtain the
following proposition.

\begin{Prop}
The controls which minimize the cost function $J$ given by $(5.4)$
and steers the system $(5.3)$ from $X(0)=X_{0} $ at $t=0$ to
$X(t_{f})= X_{f} $ at $t=t_{f}$ are
 given by $~u_{1}=\dfrac{z_{1}}{c_{1}},~
u_{2}=\dfrac{z_{2}}{c_{2}} $ where $ z_{i},~i=\overline{1,4}$ are
solutions of the system:
\begin{equation}
 \dot{z}_{1}  = \dfrac{1}{c_{2}}
z_{2} z_{3},~~~ \dot{z}_{2}  = - \dfrac{1}{c_{1}} z_{1} z_{3},~~~
\dot{z}_{3}  = - \dfrac{1}{c_{1}} z_{1} z_{4},~~~\dot{z}_{4}
=0.~~~~~\hfill\Box \label{(5.5)}
\end{equation}
\end{Prop}

It is easy to see from the equations $(5.5)$ that $z_{4}=k$ (
$k=constant$) and so the system $(5.5)$ can be written in the
equivalent form:
\begin{equation}
 \dot{z}_{1}  = \dfrac{1}{c_{2}}
z_{2} z_{3},~~~ \dot{z}_{2}  = - \dfrac{1}{c_{1}} z_{1} z_{3},~~~
\dot{z}_{3}  = - \dfrac{k}{c_{1}} z_{1}. \label{(5.6)}
\end{equation}

We observe that $(5.6)$ is a differential system which belongs to
Maxwell-Bloch top system. For the study of geometrical and
dynamical properties of the system $(5.6)$ we apply the results
given in Sections 2, 3 and 4. To this end  we consider the
following change of variables:
\begin{equation}
 z_{1} = y_{2},~~~ z_{2}= y_{3},~~~ z_{3}= y_{1}. \label{(5.7)}
\end{equation}

Using the relations $(5.7)$, the system $(5.6)$ reads:

\begin{equation}
 \dot{y}_{1}  = -\dfrac{k}{c_{1}}
y_{2},~~~ \dot{y}_{2}  = \dfrac{1}{c_{2}} y_{1} y_{3},~~~
\dot{y}_{3}  = - \dfrac{1}{c_{1}} y_{1} y_{2}. \label{(5.8)}
\end{equation}

Applying Proposition 2.1 for the system $(5.8)$ and the relation
$(5.7)$ we obtains the following proposition.

\begin{Prop}
 The functions $ \widetilde{H}, \widetilde{C} \in C^{\infty}({\bf R}^{3}, {\bf R}) $  given by:
\begin{equation}
 \widetilde{H}(z_{1},z_{2}, z_{3}) = -
\displaystyle\frac{k}{2c_{1}}\left(z_{1}^{2}+
\displaystyle\frac{c_{1}}{c_{2}} z_{2}^{2}\right)~~\hbox{and}~~
\widetilde{C}(z_{1}, z_{2}, z_{3}) = z_{2}-
\displaystyle\frac{1}{2 k} z_{3}^{2} \label{(5.9)}
\end{equation}
are constants of the motion for the dynamics (5.6).\hfill$\Box$
\end{Prop}

\begin{Rem} {\rm From Proposition $5.2$ it follows that} the trajectories of the dynamics (5.6) in
the phase space ${\bf R}^{3}$ are the intersections of the
surfaces:
\[
- \displaystyle\frac{k}{c_{1}}\left(z_{1}^{2}+
\displaystyle\frac{c_{1}}{c_{2}} z_{2}^{2}\right) = 2
\widetilde{H},~~~ 2 z_{2}- \displaystyle\frac{1}{k} z_{3}^{2}=
2\widetilde{C},
\]
where $~\widetilde{H}=constant~$ and $~\widetilde{C}=constant .$
\hfill$\Box$
\end{Rem}

Next proposition shows that the dynamics of the system $(5.6)$
reduces to pendulum dynamics.

\begin{Prop}
The solution of the system $(5.6)$, restricted to the constant
level surface defined by:
\begin{equation}
z_{1}^{2} + \displaystyle\frac{c_{1}}{c_{2}} z_{2}^{2}= 2 H =
constant, ~~~ H > 0 \label{(5.10)}
\end{equation}
is
\begin{equation}
 z_{1}(t) =
\sqrt{2H}\cdot \cos\theta(t),~~ z_{2}(t) =
\sqrt{\dfrac{c_{2}}{c_{1}}} \sqrt{2H}\cdot \sin \theta(t),~~
 z_{3}(t) =-\sqrt{c_{1}c_{2}}\cdot \dot{\theta}(t)\label{(5.11)}
\end{equation}
where $\theta(t)$ is a solution of the pendulum equation:
\begin{equation}
\ddot{\theta}(t)= \dfrac{k}{c_{1}^{2}}
\sqrt{\frac{c_{1}}{c_{2}}}\sqrt{2H} \cdot \cos
\theta(t).\label{(5.12)}
\end{equation}
\end{Prop}
{\it Proof.} We apply Proposition 2.2 for the system $(5.8)$ and
 replace $y_{i}$ with $z_{j}$ according with $(5.7).\hfill\Box$

\begin{Prop}
The system $(5.6)$ has the Hamilton-Poisson realization $({\bf
R}^{3}, \Pi, \widetilde{H})$ with the Hamiltonian $\widetilde{H}$
and Casimir  $ \widetilde{C}\in C^{\infty}({\bf R}^{3},{\bf R})$
given by $(5.9)$ and
\[
\Pi(z_{1}, z_{2} , z_{3}) = \left ( \begin{array}{ccc} 0 & -\dfrac{1}{k}z_{3} & -1\\
\\
 \dfrac{1}{k}z_{3} & 0 &  0\\ \\
1 & 0 &  0 \\
\end{array}\right ).
\]
\end{Prop}

{\it Proof.} The system $(5.8)$ is a Maxwell-Bloch top system with
$ b=( -k /c_{1}, 1 / c_{2}, - 1 / c_{1}).$

Applying Proposition 3.1 and the relation $(3.2)$ we obtain that
$(5.8)$ has the Hamilton-Poisson realization $({\bf R}^{3},
\Pi_{1}(y), H_{1}(y))$ with the Casimir $C_{1}(y)$, where
\[
\Pi_{1}(y) = \left ( \begin{array}{ccc} 0 & 1 & 0\\
-1 & 0 & - \dfrac{1}{k}y_{1}\\
0 & \dfrac{1}{k}y_{1} &  0 \\
\end{array}\right ),~H_{1}(y)=-\dfrac{k}{2 c_{1}} \left( y_{2}^{2}
+\dfrac{c_{1}}{c_{2}}y_{3}^{2}\right),~C_{1}(y)=y_{3}-
\dfrac{1}{2k} y_{1}^{2}.
\]

Replacing in $H_{1}(y)$ and $C_{1}(y)$  the variables $y_{i}$ with
$z_{j}$ by the change of variables given by $(5.7)$, we find the
Hamiltonian $\widetilde{H}(z)$ and Casimir $\widetilde{C}(z)$
given in $(5.9)$.

To determine $\Pi(z_{1}, z_{2}, z_{3})=(\{z_{i},z_{j}\})$ we use
$\Pi_{1}(y) $ and $(5.7)$. We have $\{z_{1}, z_{2}\}=\{y_{2},
y_{3}\}_{1}=-1 /k y_{1}=- 1 /k z_{3},~\{z_{1}, z_{3}\}=\{y_{2},
y_{1}\}_{1}=-1,~\{z_{2}, z_{3}\}= \{y_{3}, y_{1}\}_{1}= 0.$ Thus
we obtain the matrix $\Pi(z_{1}, z_{2}, z_{3})$.\hfill$\Box$

The equilibrium states of the dynamics $(5.6)$ are
$~\widetilde{e}_{0} = (0, 0, 0),~ \widetilde{e}_{2}^{m} = (0, m,
0)$ and $ \widetilde{e}_{3}^{m}= (0, 0, m)$ for all $~ m \in {\bf
R}^{\ast}.$

 To establish the stability of equilibrium
 states for the dynamics $(5.6)$ we apply Propositions 4.2-4.4.

\begin{Prop}
$(i)~\widetilde{e}_{0} $ and $ \widetilde{e}_{3}^{m}$ for $ m \in
{\bf R}^{\ast} $ are nonlinear stable.

$(ii)~\widetilde{e}_{2}^{m}$ for $ m \in {\bf R}^{\ast} $ are
nonlinear stable if $~k m > 0~ $ and unstable if $~k m < 0
.$\hfill$\Box$
\end{Prop}

{\bf Conclusions}. In this paper we have presented the geometric
and dynamical properties of the Maxwell-Bloch top system $(2.1).$
The denomination used is justified by the fact that the $3D$
 Maxwell-Bloch equations belongs to respective
family.\hfill$\Box$

{\bf Acknowledgments.} The author has very grateful to be
reviewers for their comments and suggestions.\\[-0.3cm]

Author's adress\\

Mihai Ivan\\[0.1cm]
West University of Timi\c soara\\
Seminarul de Geometrie \c si Topologie. Department of Mathematics\\
4, B-dul V. P{\^a}rvan, 300223, Timi\c{s}oara, Romania\\
E-mail: ivangm31@yahoo.com\\

\end{document}